\documentclass[a4paper,12pt]{amsart}

\usepackage{amscd}
\usepackage{amsfonts}
\usepackage{amsmath}
\usepackage{amssymb}
\usepackage{latexsym}
\usepackage{amsthm,color}
\usepackage[all]{xy}

\usepackage{epsfig}
\usepackage{graphicx}
\usepackage{color}
\usepackage{tikz}
\usetikzlibrary{calc}
\newtheorem{lemma}{Lemma}[section]

\newtheorem{theorem}{Theorem}
\newtheorem{corollary}{Corollary}

\theoremstyle{definition}
\newtheorem{definition}{Definition}
\theoremstyle{definition}

\theoremstyle{definition}

\newtheorem{note}{Note}[section]
\newtheorem{property}{Property}[section]

\newcommand{\p}{\mathbb{P}}

\begin{document}
\title[
Spherical orthotomic curve-germs]
{Spherical orthotomic curve-germs}
\author[X.~Liu]{Xihe Liu}
\address{
Graduate School of Environment and Information Sciences, 
Yokohama National University, {Yokohama 240-8501,} Japan
}
\email{liu-xihe-tx@ynu.jp}
\author[T.~Nishimura]{Takashi Nishimura
}
\address{
Research Institute of Environment and Information Sciences,  
Yokohama National University, 
Yokohama 240-8501, Japan}
\email{nishimura-takashi-yx@ynu.jp}
\begin{abstract}
In this paper, it is shown that 
for an $n$-dimensional spherical unit speed curve    
$\gamma: I\to S^n$, a given point $P \in S^n$ and a point 
$s_0$ of the open 
interval $I$, 
the spherical orthotomic curve-germ 
$ort_{\gamma, P}: (I, s_0)\to S^n$
of $\gamma$ relative to $P$ is 
$\mathcal{L}$-equivalent to the spherical pedal curve-germ   
$ped_{\gamma, P}: (I, s_0)\to S^n$
of $\gamma$ relative to $P$  (resp., the spherical dual curve-germ  
${\bf u}_n: (I, s_0)\to S^n$ of $\gamma$) 
if and only if $P\ne \pm{\bf u}_n(s_0)$ (resp., if $P= \pm{\bf u}_n(s_0)$).    
\end{abstract}


\date{}

\maketitle

\maketitle

\section{Introduction\label{section 1}}
Throughout this paper, let $I$, $S^n$ be an open interval  
and the unit sphere in 
$\mathbb{R}^{n+1}$ $(n\ge 2)$ respectively.      
We are concerned with singularities of spherical orthotomic curve.    
In the case of $\mathbb{R}^2$, plane orthotomic curves 
are classically well-studied 
in differential geometry  and have interesting applications 
(for instance, see \cite{brucegiblin}).    
Given a plane unit-speed curve $\gamma: I\to \mathbb{R}^2$ and a point 
$P$ of $\mathbb{R}^2$,  
the pedal curve of $\gamma$ relative to $P$ and 
the orthotomic curve of $\gamma$ relative to $P$ 
are  naturally defined; and it is easily seen that the orthotomic curve is 
$\mathcal{L}$-equivalent to the pedal curve.    
Thus, in $\mathbb{R}^2$, the study of singularities of the orthotomic curve 
comes down to the study of singularities of the pedal curve 
(for details, see Section \ref{section 5}).    
Since the recognition problem for generic singularities 
of spherical pedal curve have been already solved 
(see \cite{geomdedicata, demonstratio}), in the case of $S^n$ as well,  
it is expected to have a similar relationship as the plane curve case.   
However, in the case of $S^n$, 
for a given spherical unit-speed curve $\gamma: I\to S^n$ and a given 
point $P$ of $S^n$,  
the spherical pedal curve of $\gamma$ relative to $P$ is not always defined 
although the spherical orthotomic curve of $\gamma$ relative to $P$ 
is always defined.     
Moreover, even if the spherical pedal curve is defined, it is 
uncertain about the $\mathcal{L}$-equivalence between the spherical 
orthotomic curve and  the spherical pedal curve.      
Therefore, 
generic singularities of spherical orthotomic curve is wrapped in mystery.      
\par 
The purpose of this paper is to clarify the mystery surrounding spherical 
orthotomic curve-germs.   
\begin{definition}\label{def1}(\cite{geomdedicata})
{\rm
A regular  curve $\gamma :I \rightarrow S^n$ is called 
a $\mathit{spherical}$ $\mathit{unit}$ $\mathit{speed}$ $\mathit{curve}$ 
if the following $\mathbf{u}_i:I \rightarrow S^n$ ($1\leqslant i \leqslant n-1$) 
is inductively well-defined:
 \begin{gather*}
  \mathbf{u}_{-1}(s) = \mathbf{0}, \quad \mathbf{u}_{0}(s) = \gamma(s), 
\quad \parallel\mathbf{u}_{0}'(s)\parallel \equiv 1 , \quad\kappa_0(s)\equiv 0,\\
  \mathbf{u}_{i}(s)=\frac{\mathbf{u}_{i-1}'(s)+\kappa_{i-1}(s)\mathbf{u}_{i-2}(s)}
{\parallel\mathbf{u}_{i-1}'(s)+\kappa_{i-1}(s)\mathbf{u}_{i-2}(s)\parallel}\quad 
(1\leqslant i \leqslant n-1)\\
  \kappa_i(s) =\parallel\mathbf{u}_{i-1}'(s)+
\kappa_{i-1}(s)\mathbf{u}_{i-2}(s)\parallel 
\quad (1\leqslant i \leqslant n-1)
\end{gather*}
}
\end{definition}
Notice that the inductive assumption $\kappa_i(s)>0$ for any 
$i$ $(1\leqslant i \leqslant n-1)$ and any $s\in I$ 
are not so strong (see \cite{geomdedicata}).   
In \cite{geomdedicata}, the following has been shown .
\begin{lemma}\label{lemma1}\rm{(\cite{geomdedicata})}
For any $s$ $\in$ $I$ and any integers $i$,$k$ such that $-1\leqslant i < k 
\leqslant n-1$, the following hold:
\begin{enumerate}
\item[(1)] $\mathbf{u}_i(s)\boldsymbol{\cdot}\mathbf{u}_k(s) = 0$, 
\item[(2)] $\mathbf{u}_i(s)\boldsymbol{\cdot}\mathbf{u}_k'(s) = 0\quad(i<k-1)$,
\item[(3)] $\mathbf{u}_{k-1}(s)\boldsymbol{\cdot}\mathbf{u}_k'(s) 
= -\kappa_{k}(s)$.
\end{enumerate}
Here, the dot in the center stands for the scalar product of two vectors of $
\mathbb{R}^{n+1}$.
\end{lemma}
\begin{definition}\label{def2}
(\cite{geomdedicata}, see also \cite{arnold, porteousbanach, porteous})
{\rm
Let $\gamma :I \rightarrow S^n$ be a spherical unit speed curve.
\begin{enumerate}
\item[(1)]  The mapping $\mathbf{u}_{n}:I \rightarrow S^n$, 
called the {\it spherical dual curve} of $\gamma$, is defined by
\[
\det\left(\mathbf{u}_{0}(s),\dots,\mathbf{u}_{n-1}(s),\mathbf{u}_{n}(s)\right)=1.
\]
\item[(2)]    
The function $\kappa_n: I \to \mathbb{R}$ is defined by
\[
\kappa_n(s)=\mathbf{u}_{n-1}'(s)\cdot\mathbf{u}_{n}(s).
\]
\end{enumerate}
}
\end{definition}
By Definition \ref{def2}, it is clear that 
the set consisting of $(n+1)$-vectors 
${\bf u}_0(s), \ldots, {\bf u}_n(s)$ 
is a unit orthogonal basis of $\mathbb{R}^{n+1}$.  
The spherical pedal curve and 
the spherical orthotomic curve are defined as follows:
\begin{definition}[\cite{geomdedicata}]\label{def3}
{\rm
Let $\gamma :I \rightarrow S^n$ be 
a spherical unit speed curve.    
Moreover, let 
\[
P=\sum_{i=0}^n \left(P\cdot {\bf u}_i(s)\right){\bf u}_i(s) 
\] 
be a point of $S^n$ such that $\left(P\cdot {\bf u}_n(s)\right)\ne \pm 1$ 
for any $s\in I$.  
Then, the {\it spherical pedal curve} of $\gamma$ relative to $P$, 
denoted by $ped_{\gamma,P}:I \rightarrow S^n$ is defined as follows:   
\[
ped_{\gamma,P}(s)=
\frac{\sum_{i=0}^{n-1}\left(P\cdot {\bf u}_i(s)\right){\bf u}_i(s)}
{\left\|\sum_{i=0}^{n-1}\left(P\cdot {\bf u}_i(s)\right){\bf u}_i(s)\right\|}.    
\]
}
\end{definition}
Notice that for a point $P$ such that there exists $s_0\in I$ 
satisfying $\left(P\cdot {\bf u}_n(s_0)\right)=1$, the spherical pedal curve  
$ped_{\gamma,P}$ is not well-defined when $s=s_0$.
By Definition \ref{def3}, it is easily seen that the following holds:   
\begin{lemma}[Lemma 3.1 of \cite{geomdedicata}]\label{pedal}
\[
ped_{\gamma,P}(s)=
\frac{1}{\sqrt{1-\left(P\cdot {\bf u}_n(s)\right)^2}}
\left(P-\left(P\cdot {\bf u}_n(s)\right){\bf u}_n(s)\right).   
\]
\end{lemma} 
\begin{definition}\label{def4}
{\rm
Let $\gamma :I \rightarrow S^n$ be 
a spherical unit speed curve.    
Moreover, let 
\[
P=\sum_{i=0}^n \left(P\cdot {\bf u}_i(s)\right){\bf u}_i(s) 
\] 
be a point of $S^n$.  
Then, the {\it spherical orthotomic curve} of $\gamma$ relative to $P$, 
denoted by $ort_{\gamma,P}:I \rightarrow S^n$ is defined as follows:   
\[
ort_{\gamma,P}(s)=\sum_{i=0}^{n-1}\left(P\cdot {\bf u}_i(s)\right){\bf u}_i(s) 
-\left(P\cdot {\bf u}_n(s)\right){\bf u}_n(s).   
\]
}
\end{definition}
Notice that for a given unit-speed curve $\gamma$ and a point $P\in S^n$, 
the spherical orthotomic curve of $\gamma$ relative to $P$ 
is always defined while the spherical pedal curve for $\gamma$ and $P$ 
is not always defined.   
Notice also that when $P= \pm\mathbf{u}_n(s_0)$ , it follows that 
$ort_{\gamma,P}(s_0)=- P$.
\begin{definition}\label{def4}
{\rm
Two curve-germs $f,g:(I, 0)\to S^n$ are said to be 
{\it $\mathcal{L}$-equivalent}  
if there exists a germ of $C^\infty$ diffeomorphism 
$\psi:(S^n, f(0))\to(S^n, g(0))$ such that the equality $g=\psi\circ f$ holds.
}
\end{definition}
The main result of this paper is the following:   
\begin{theorem}\label{thm1}
Let $\gamma: I\to S^n$ be a spherical unit speed curve and let 
$P$ be a point of $S^n$.    
 Then, for any $s_0\in I$, the following hold.   
\begin{enumerate}
\item[(1)]   
The orthotomic curve-germ $ort_{\gamma, P}: (I, s_0)\to S^n$ is
$\mathcal{L}$-equivalent to
 $ped_{\gamma, P}: (I, s_0)\to S^n$ if and only if 
$P\ne \pm{\bf u}_n(s_0)$.
\item[(2)]   
The orthootomic curve-germ 
$ort_{\gamma, P}: (I, s_0)\to S^n$ is
$\mathcal{L}$-equivalent to
 ${\bf u}_n: (I, s_0)\to S^n$ if $P=\pm{\bf u}_n(s_0)$.
\end{enumerate}
\end{theorem} 
Notice the following: 
\begin{note}\label{note1}
\begin{enumerate}
\item[(1)] The pedal curve-germ $ped_{\gamma, P}: (I, s_0)\to S^n$ 
is not well-defined in the case $P=\pm{\bf u}_n(s_0)$.
\item[(2)] The converse of the assertion (2) of Theorem \ref{thm1} 
does not hold in general 
because the curve-germ $ped_{\gamma, P}: (I, s_0)\to S^n$ 
is generically $\mathcal{L}$-equivalent to the dual curve-germ 
${\bf u}_n: (I, s_0)\to S^n$ even if $P\ne \pm{\bf u}_n(s_0)$ holds.   
\item[(3)] Even if $P\ne \pm{\bf u}_n(s_0)$ holds, there are examples 
that $ped_{\gamma, P}: (I, s_0)\to S^n$ 
is not $\mathcal{L}$-equivalent to ${\bf u}_n: (I, s_0)\to S^n$ 
(see Corollary \ref{coro1-1} in Section \ref{section 2}).   
\end{enumerate}
\end{note}
\par  
Since singularities of ${\bf u}_n: (I, s_0)\to S^n$ and 
singularities of $ped_{\gamma, P}$ have been relatively-well investigated, 
thanks to Theorem \ref{thm1}, normal forms of 
generic singularities of $ort_{\gamma, P}$ can be obtained systematically and 
recognizably (for details, see 
Section \ref{section 2}).   
\par 
\bigskip 
This paper is organized as follows. In Section \ref{section 2}, 
applications of Theorem \ref{thm1} are given.    
In Section \ref{section 3},  
preliminaries for the proof of Theorem 1 are given.  
Theorem \ref{thm1} is proved in Section \ref{section 4}.
Finally, for readers' convenience, 
plane pedal curves and plane orthotomic curves are 
reviewed in Section \ref{section 5}.    
\section{Applications of Theorem \ref{thm1}}\label{section 2}
\subsection{Application 1}
Singularities of spherical orthotomic curve-germs have been less-studied.   
To the best of authors' knowledge, only \cite{xiong} is the literature 
on this topic.   
Although \cite{xiong} is the work of pioneer and 
his results are interesting, 
the tool used there is the standard one 
in Singularity Theory given by \cite{brucegiblin}, 
and thus the results obtained in \cite{xiong} are  
just the $\mathcal{V}$-equivalence of images of parametrization-germ, 
not the $\mathcal{L}$-equivalence of parametrization-germs.    
It is significant to strengthen the results of \cite{xiong} to results on    
$\mathcal{L}$-equivalence of parametrization-germs,  
and also to investigate the case that 
the given point $P$ is over the curve $\gamma$ (this case 
has not been treated in \cite{xiong}).   
Moreover, we wanted to generalize his results in higher dimensions.   
By Theorem \ref{thm1}, all of the above purposes can be realized.   
\subsection{Application 2}
\begin{definition}\label{def6}
{\rm
For any $i$ $(-1\leqslant i\leqslant n)$, define $S^{i}_{\mathbf{u}_i(s)}$ by
\[
S^{i}_{\mathbf{u}_i(s)} = (S^n - \{ \pm \mathbf{u}_{n}(s) \}) 
\cap \mathop{\mathord{\sum}}\limits_{j=-1}^{i}\mathbb{R}\mathbf{u}_j(s).
\]
}
\end{definition}
Combining Theorem 1.1 of \cite{geomdedicata} 
and Theorem \ref{thm1}, the following corollary is obtained.
\begin{corollary}\label{coro1}
Suppose that the dual curve ${\bf u}_n$ is non-singular and 
$P$ is located outside ${\bf u}_n(I)$.   Then, 
singularity types of $ort_{\gamma, P}$ depend only on the location of $P$.   
Thus, for instance, we have the following:
\begin{enumerate}
\item[(1)] The point $P$ is inside 
$S_{\mathbf{u}_{n}(s_0)}^{n}-S_{\mathbf{u}_{n-2}(s_0)}^{n-2}$ 
if and only if the map-germ $ort_{\gamma,P}:(I,s_0) \rightarrow S^n$ 
is $\mathcal{L}$-equivalent to the map-germ given by 
$s \longmapsto (s,0,\ldots,0).$
\item[(2)] For any i $(0\leqslant i  \leqslant n-2)$, 
the point $P$ is inside $S_{\mathbf{u}_{i}(s_0)}^{i}-S_{\mathbf{u}_{i-1}(s_0)}^{i-1}$ 
if and only if the map-germ $ort_{\gamma,P}:(I,s_0) \rightarrow S^n$ 
is $\mathcal{L}$-equivalent to the map-germ given by the following:
    \[
    s \longmapsto (\underbrace{s^{n-i},s^{n-i+1},\ldots,s^{2n-2i-1}}_{(n-i)
\hspace{1mm}elements},\underbrace{0,\ldots,0}_{i\hspace{1mm}elements}).
    \]
\end{enumerate}
\end{corollary}
Notice that in the assertion (2) of Corollary \ref{coro1}, 
the orthotomic curve-germ 
$ort_{\gamma,P}:(I,s_0) \rightarrow S^n$ must be singular while the 
dual curve-germ ${\bf u}_n: (I, s_0)\to S^n$ is non-singular.  
Therefore, as a by-product of Corollary \ref{coro1}, we have the following 
Corollary \ref{coro1-1}.   
\begin{definition}\label{def5}
{\rm
Two curve-germs $f,g:(I, s_0)\to S^n$ are said to be 
{\it $\mathcal{A}$-equivalent}  
if there exist germs of $C^\infty$ diffeomorphism 
$\phi:(I,s_0)\to(I,s_0)$ and $\psi:(S^n,f(s_0))\to (S^n, g(s_0))$ 
such that the equality $g\circ \phi=\psi\circ f$ holds.    
}
\end{definition}
\begin{corollary}\label{coro1-1}
Suppose that the dual curve ${\bf u}_n$ is non-singular and 
the point $P$ is located inside 
$S_{\mathbf{u}_{i}(s_0)}^{n-2}$.   Then, 
the following hold: 
\begin{enumerate}
\item[(1)]   The orthotomic curve-germ 
$ort_{\gamma, P}: (I, s_0)\to S^n$ 
is $\mathcal{L}$-equivalent to 
$ped_{\gamma, P}: (I, s_0)\to S^n$. 
\item[(2)]   The orthotomic curve-germ 
$ort_{\gamma, P}: (I, s_0)\to S^n$ is never $\mathcal{A}$-equivalent to 
${\bf u}_n: (I, s_0)\to S^n$.     
\end{enumerate}
\end{corollary}
\par 
\smallskip 
Combining Theorems 2 and 3 of \cite{demonstratio} 
and the above Theorem \ref{thm1}, the following corollaries are obtained.
\begin{definition}\label{def7}
{\rm
The function $f: I\to \mathbb{R}$ 
is said to {\it have $\mathit{A_k}$-$\mathit{type}$ 
$\mathit{singularity}$} ($k \geqslant 0$) at $s_0\in I$ 
if $f(s_0)=f'(s_0)=\ldots=f^{(k)}(s_0)=0$ and $f^{(k+1)}(s_0)\ne 0$. 
Thus, $\mathit{A_0}$-$\mathit{type}$ $\mathit{singularity}$ 
means $f(s_0)=0$ $f'(s_0)\ne 0$.
}
\end{definition}
\begin{corollary}\label{coro2}
Let $\gamma:I \rightarrow S^n$ be a spherical unit speed curve. 
Suppose that for $s_0\in I$ 
$P$ is located inside $S^{n}_{\mathbf{u}_{n}(s_0)}-
S^{n-1}_{\mathbf{u}_{n-1}(s_0)}$. 
Then the following hold.
\begin{enumerate}
\item[(1)] If $\kappa_n$ has an $A_k$-type singularity at $s_0$ 
$(0\leqslant k \leqslant n-2)$, 
then the orthotomic curve-germ $ort_{\gamma,P}:(I,s_0) 
\rightarrow S^n$ is $\mathcal{L}$-equivalent 
to the map-germ given by
  \[
    s \longmapsto (\underbrace{s^{k+2},s^{k+3},\ldots,s^{2k+3}}_{(k+2)
\hspace{1mm}
elements},\underbrace{0,\ldots,0}_{n-k-2\hspace{1mm}elements}).
    \]
\item[(2)] If $\kappa_n$ has an $A_{n-1}$-type singularity at $s_0$, 
then the orthotomic curve-germ 
$ort_{\gamma,P}:(I,s_0) \rightarrow S^n$ is 
$\mathcal{A}$-equivalent to the 
map-germ given by
    \[
    s \longmapsto (s^{n+1},s^{k+2},\ldots,s^{2n}).
    \]
\end{enumerate}
\end{corollary}
\begin{corollary}\label{coro3}
Let $\gamma:I\rightarrow S^n$ be a spherical unit speed curve. 
Suppose that $\kappa_n$ has an $A_0$-type singularity at $s_0$.    
Then the following 
hold.
\begin{enumerate}
\item[(1)]The pedal point $P$ is inside 
$S^{n}_{\mathbf{u}_{n}(s_0)}-S^{n-1}_{\mathbf{u}_{n-1}(s_0)}$ 
if and only if the orthotomic curve-germ 
$ort_{\gamma,P}:(I,s_0) \rightarrow S^n$ is 
$\mathcal{L}$-equivalent to the map-
germ given by
    \[
    s \longmapsto (s^{2},s^{3},\ldots,0).
    \]
\item[(2)] For any $i$ $(1 \leqslant i \leqslant n-1)$, the point $P$ is 
inside $S^{n-i}_{\mathbf{u}_{n-i}(s_0)}-S^{n-i-1}_{\mathbf{u}_{n-i-1}(s_0)}$ 
if and only if 
the orthotomic curve-germ $ort_{\gamma,P}:(I,s_0) \rightarrow S^n$ is 
$\mathcal{A}$-equivalent to the map-germ given by
    \[
    s \longmapsto (\underbrace{s^{i+1},s^{i+3},s^{i+4}\ldots,s^{2i+3}}_{(i-1)
\hspace{1mm}elements},\underbrace{0,\ldots,0}_{(n-i-1)\hspace{1mm}elements}).
    \]
\item[(3)] The point $P$ is inside $S^{0}_{\mathbf{u}_{0}(s_0)}-
S^{-1}_{\mathbf{u}_{-1}(s_0)}$ if and only if the orthotomic curve-germ 
$ort_{\gamma,P}:(I,s_0) \rightarrow S^n$ is 
$\mathcal{A}$-equivalent to the map-germ given by
     \[
    s \longmapsto (s^{n+1},\underbrace{s^{n+3},s^{n+4}\ldots,s^{2n+1}}_{(n-1)
\hspace{1mm}elements}).
    \]
\end{enumerate}
\end{corollary}
\section{Preliminaries\label{section 3}}
\begin{figure}[h]
\centering
\begin{tikzpicture}
  [
    scale=1.58,
    >=stealth,
    point/.style = {draw, circle,  fill = black, inner sep = 1pt},
    dot/.style   = {draw, circle,  fill = black, inner sep = .2pt},
  ]

  \def\rad{1}
  \node (origin) at (0,0) [point, label = {below:$O$}]{};
  \draw (origin) circle (\rad);

  \node (n1) at +(50:\rad) [point, label = above:\qquad$ped_{\gamma,P}(s)$] {};
  \node (n2) at +(100:\rad) [point, label = above:$P$] {};
  \node (n3) at +(0:\rad) [point, label = {below right:$ort_{\gamma,P}(s)$}] {};

  \draw (n3) -- node (a) {} (n1);
  \draw (n3) -- node (b) {} (n2);
  \draw (n2) -- (n1);
  \draw (origin) -- (n1);
  \draw (origin) -- (n2);
  \draw (origin) -- (n3);
 \def\ralen{.5ex}  
  \foreach \inter/\first/\last in {b/n2/origin}
    {
      \draw let \p1 = ($(\inter)!\ralen!(\first)$), 
                \p2 = ($(\inter)!\ralen!(\last)$),  
                \p3 = ($(\p1)+(\p2)-(\inter)$)      
            in
              (\p1) -- (\p3) -- (\p2)               
              ($(\inter)!.5!(\p3)$) node [dot] {};  
   }
   \draw(0.096418141452981,0.114906666467847)arc (50:100:0.15) node[above]
{\quad$\theta$};
\end{tikzpicture}
\caption{
$ped_{\gamma, P}(s)$ and 
$ort_{\gamma, P}(s)$ in $S^n$.}   
\label{figure2}
\end{figure}
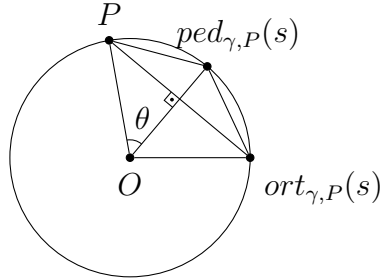

\begin{lemma}\label{lemma3}\rm{(\cite{geomdedicata})}
$ped'_{\gamma,P}(s) = \mathbf{0} \Leftrightarrow$
$\kappa(s) = 0$ or $P \in$ $S_{\mathbf{u}_{n-2}(s)}^{n-2}$ = $(S^n - \{ \pm 
\mathbf{u}_{n}(s) \})$ $\cap$ $\mathop{\mathord{\sum}}\limits_{j=-1}^{n-2}
\mathbb{R}\mathbf{u}_j(s)$.
\end{lemma}
\begin{lemma}\label{lemma4}\rm{(\cite{geomdedicata})}
$\mathbf{u}'_n(s) = -\kappa_n(s)\mathbf{u}_{n-1}(s)$.
\end{lemma}
From Lemma \ref{lemma4}, it follows that Lemma \ref{lemma3} is equivalent to 
the following:
\begin{lemma}\label{lemma5}\rm{(\cite{geomdedicata})}
$ped'_{\gamma,P}(s) = \mathbf{0} \Leftrightarrow$
$\mathbf{u}'_n(s) = \mathbf{0}$ or $P \in$ $S_{\mathbf{u}_{n-2}(s)}^{n-2}$.
\end{lemma}
\begin{definition}\label{def9}(\cite{geomdedicata})
{\rm
Let $\gamma :I \rightarrow S^n$ be a $\mathit{spherical}$ $\mathit{unit}$ $
\mathit{speed}$ $\mathit{curve}$. Let $P$ be a point inside 
$S^n - \{ \pm \mathbf{u}_{n}(s) \vert s 
\in I \}$. Then, $\Psi_P:S^n-\{\pm P\}$ $\rightarrow$ $S^n$ is defined by:
\[
\Psi_P(\mathbf{x})=\frac{1}{\sqrt{1-(P\boldsymbol{\cdot}\mathbf{x})^2}}(P-(P
\boldsymbol{\cdot}\mathbf{x})\mathbf{x}).
\]
}
\end{definition}
From Definition \ref{def9} and Lemma \ref{pedal}, it follows that  
$ped_{\gamma,P}$ is factored into two mappings as follows:
\begin{property}(\cite{geomdedicata})\label{property1}
\[
ped_{\gamma,P}(s)=\Psi_P \circ \mathbf{u}_{n}(s).
\]
\end{property}
By Property \ref{property1}, Lemma \ref{lemma3} 
is equivalent to the following Lemma \ref{lemma6}:
\begin{lemma}\label{lemma6}\rm{(\cite{geomdedicata})}
$ped'_{\gamma,P}(s)=\mathbf{0}$ $\Leftrightarrow$ $\mathbf{u}'_{n}(s)=\mathbf{0}$ or 
$\mathbf{u}_{n}(s)$ is a singular point of $\Psi_P$.
\end{lemma}
Figure \ref{figure2} shows the relationship between spherical pedal 
curves and spherical orthotomic curves when spherical pedal curves exist.
By Figure \ref{figure2}, the following equality holds.    
Notice that $\cos \theta$ is more than $0$ and 
less than or equal to $1$ for the angle $\theta$ given in 
Figure \ref{figure2}.   
\[
 \frac{ort_{\gamma,P}(s)+P}{2}=(\cos \theta) ped_{\gamma,P}(s). 
\]
It is clear that 
\[
\cos\theta=P\boldsymbol{\cdot}ped_{\gamma,P}(s).
\]
Thus, we have the following.
\begin{lemma}\label{lemma7}
\[
ort_{\gamma,P}(s)=2(P\boldsymbol{\cdot}ped_{\gamma,P}(s))ped_{\gamma,P}(s)-P.
\]
\end{lemma}
\begin{definition}\label{def10}
{\rm
For any $P \in S^n$, $\Phi_P:\mathbb{R}^{n+1}\rightarrow\mathbb{R}^{n+1}$ is 
defined as follows.
\[
\Phi_P(\mathbf{x})=2(P\boldsymbol{\cdot}\mathbf{x})\mathbf{x}-P.
\]
}
\end{definition}
\begin{lemma}\label{lemma8}
Let $\gamma :I \rightarrow S^n$ be a spherical unit speed curve. Let $P$ be a 
point inside $S^n - \{ \pm \mathbf{u}_{n}(s) \vert s \in I \}$. Then, the following 
holds.
\begin{equation}
\notag
ort_{\gamma,P}(s)= \Phi_P \circ ped_{\gamma,P}(s)
= \Phi_P \circ \Psi_P \circ \mathbf{u}_{n}(s).
\end{equation}
\end{lemma}
\section{Proof of Theorem 1\label{section 4}}
\subsection{Proof of the assertion (1) of Theorem \ref{thm1}}
Since the pedal curve-germ $ped_{\gamma, P}: (I, s_0)\to S^n$ can be defined 
only when $P\ne \pm{\bf u}_n(s_0)$, it is sufficient to show that 
$ort_{\gamma, P}$ is $\mathcal{L}$-equivalent to 
$ped_{\gamma, P}$ if $P\ne \pm{\bf u}_n(s_0)$.   
\par 
From Lemma \ref{lemma8}, the following Lemma \ref{lemma9} holds:   .
\begin{lemma}\label{lemma9}
$ort'_{\gamma,P}(s)=\mathbf{0}$ 
$\Leftrightarrow$ $ped'_{\gamma,P}(s)=\mathbf{0}$ 
or $ped_{\gamma,P}(s)$ is a singular point of $\Phi_P$.
\end{lemma}
\begin{lemma}\label{lemma10}
The following three hold.
\begin{enumerate}
\item[(1)] $\Phi_P(S^n)$=$S^n$,
\item[(2)] $\mathbf{x} \in \mathbb{R}^{n+1}$ is a singular point of 
$\Phi_P$ $\Leftrightarrow$ $P\boldsymbol{\cdot}\mathbf{x}= 0 $,
\item[(3)] $\mathbf{x} \in S^n$ is a singular point of $\Phi_P \vert S^n$ 
$\Leftrightarrow$ $n \geqslant 2$ and $P\boldsymbol{\cdot}\mathbf{x}= 0$.
\end{enumerate}
\end{lemma}
$Proof$ $of$ $the$ $assertion$ $(1)$ $of$ $Lemma$ $\ref{lemma10}.$
\quad 
In order to show $\Phi_P(S^n)\subset S^n$, 
take an arbitrary point $\mathbf{x}$ of $S^n$. Then, 
\begin{eqnarray*}
\Phi_P(\mathbf{x}) \boldsymbol{\cdot} \Phi_P(\mathbf{x}) 
& = & 
(2(P \boldsymbol{\cdot}\mathbf{x})\mathbf{x}-P)\boldsymbol{\cdot}(2(P 
\boldsymbol{\cdot}\mathbf{x})\mathbf{x}-P) \\ 
{ } &= & 
4((P \boldsymbol{\cdot}\mathbf{x})^2 \parallel \mathbf{x} \parallel^2)-4P
\boldsymbol{\cdot}((P \boldsymbol{\cdot}\mathbf{x})\mathbf{x})+1\\
{ } &=& 
4(P \boldsymbol{\cdot}\mathbf{x})^2
-4(P \boldsymbol{\cdot}\mathbf{x})^2+1=1.    
\end{eqnarray*}
Therefore, $\Phi_P(S^n)$ $\subset$$S^n$.
\par 
In order to show $\Phi_P(S^n)$ $\supset$$S^n$, 
take an arbitrary point $\mathbf{y}$ 
of $S^n$.
Suppose that $\mathbf{y}\ne -P$.
Set
\[
\mathbf{x}=\frac{\frac{\mathbf{y}+P}{2}}{||\frac{\mathbf{y}+P}{2}||}.
\]
Then, it follows
\begin{eqnarray*}
2(\mathbf{x}\cdot P)\mathbf{x}-P
& = &
2\left(\frac{\frac{\mathbf{y}+P}{2}}{||\frac{\mathbf{y}+P}{2}||}\cdot P\right)
\frac{\frac{\mathbf{y}+P}{2}}{||\frac{\mathbf{y}+P}{2}||}-P \\
{ } & = &
\frac{2}{|| \mathbf{y}+P||^2}\left((\mathbf{y}\cdot P)+1\right)(\mathbf{y}+P)-P \\
{ } & = &
\frac{1}{\left(1+(\mathbf{y}\cdot P)\right)}
\left((\mathbf{y}\cdot P)+1\right)(\mathbf{y}+P)-P \\
{ } & = &
(\mathbf{y}+P)-P=\mathbf{y}.
\end{eqnarray*}
Next, suppose that $\mathbf{y}=-P$.
Let $\mathbf{x}$ be a point of $S^n$ such that $\mathbf{x}\cdot P=0$.
Then, $2(\mathbf{x}\cdot P)\mathbf{x}-P=-P=\mathbf{y}$.
Therefore, $\Phi_P(S^n)$ $\supset$$S^n$.
\hfill $\Box$
\par 
\smallskip 
$Proof$ $of$ $the$ $assertion$ $(2)$ $of$ $Lemma$ $\ref{lemma10}.$
\quad 
Set $\mathbf{e}_1=(1,0,\ldots,0)$. Let $R$ be a rotation matrix 
such that $PR=\mathbf{e}_1$. 
Let $\mathbf{x}$ be a point of $\mathbb{R}^{n+1}$.
Set $\mathbf{y}=\mathbf{x}R$. Then we have the following:
\begin{eqnarray*}
2(\mathbf{x}\cdot P)\mathbf{x}-P
& = &
2(\mathbf{y}R^{t}R\mathbf{e}^{t}_{1})\mathbf{y}R^{t}-\mathbf{e}_{1}R^{t} \\
{ } & = &
(2(\mathbf{y}\mathbf{e}^{t}_{1})\mathbf{y}-\mathbf{e}_{1})R^{t}.
\end{eqnarray*}
\noindent 
Here, $R^{t}$ (resp., $\mathbf{e}^{t}_{1}$) stands for 
the transposed matrix of $R$ (resp., $\mathbf{e}_{1}$).    
Therefore the following diagram is commutative.
\[
\begin{CD}
\mathbb{R}^{n+1} @> {\Phi_P} >> \mathbb{R}^{n+1} \\ 
@V{R}VV                                  @VV{R}V \\ 
\mathbb{R}^{n+1} @> {\Phi_{{\bf e}_1}} >> \mathbb{R}^{n+1}.  
\end{CD}
\]
\noindent 
Here, $R$ means the rotation defined by 
$\mathbf{x} \longmapsto \mathbf{x}R$. By this commutative diagram, 
it follows that $\mathbf{x}$ is a singular point of $\Phi_P$ if and only if 
$\mathbf{x}R=\mathbf{y}$ is a singular point of $\Phi_{\mathbf{e}_1}$.   
It is easily confirmed that $\mathbf{y}$ is a singular point of 
$\Phi_{\mathbf{e}_1}$ if and only if $\mathbf{e}_1\cdot\mathbf{y}=0$.  
Moreover, the following equality is easily seen.     
\[
\mathbf{e}_1\cdot\mathbf{y}=\mathbf{e}_1\mathbf{y}^t 
= \mathbf{e}_1RR^t\mathbf{y}^t 
= \left(\mathbf{e}_1R\right)\left(\mathbf{y}R\right)^t
=P\mathbf{x}^t=P\cdot \mathbf{x}.  
\]
Therefore, the proof of the assertion (2) of Lemma \ref{lemma10} completes.  
\hfill 
$\Box$ 
%
\par 
\smallskip 
$Proof$ $of$ $the$ $assertion$ $(3)$ $of$ $Lemmma$ $\ref{lemma10}.$\quad 
When n=1, 
we may set $\mathbf{x}=(x_1,x_2)=(\cos\theta, \sin\theta)$.
Then, we have the following:   
\[
\Phi_P(\mathbf{x}) = 2(P \boldsymbol{\cdot}\mathbf{x})\mathbf{x}-P
=  2\cos\theta(\cos\theta,\sin\theta)- (1,0)
= (\cos2\theta,\sin2\theta).
\]
It follows that $\Phi_P \vert_{S^1}$ is non-singular.  
\par 
Next, suppose that $n\ge 2$.   
Set $E_P=\{\mathbf{x}\in S^n\; |\; P\cdot \mathbf{x}=0\}$.    
Then, $E_P$ is an $(n-1)$-dimensional submanifold of $S^n$.   
By the assertion (2) of Lemma \ref{lemma10}, 
it follows that a point $\mathbf{x}\in S^n$ belongs to 
$E_P$ if and only if it is a singular point of $\Phi_P$.    
Notice that $\Phi_P(E_P)=-P$ 
which implies that 
for any point $\mathbf{x}$ of $E_P$, the following 
holds:   
\[
d\left(\Phi_P\right)_{\mathbf{x}}\left(T_{\mathbf{x}}E_P\right)
=\mathbf{0}. 
\]
Notice also that $\dim E_P\ge 1$ in this case.   
Therefore, a point $\mathbf{x}\in S^n$ belongs to 
$E_P$ if and only if it is a singular point of ${\Phi_P}|_{S^n}$. 
%
\hfill $\Box$
\par 
\smallskip 
For any $P\in S^n$, denote the open hemisphere centered at $P$ 
by $H(P)$.    Namely, ${H}(P)$ is defined as follows:
\[
{H}(P)=\left\{ \mathbf{x} \in S^n \vert P\boldsymbol{\cdot}
\mathbf{x} >0\right\}.
\]
Then, notice than any point of $H(P)$ is a regular point of $\Phi_P|_{S^n}$ 
by the assertion (3) of Lemma \ref{lemma10}.    Moreover, it is easily seen 
that $\Phi_P|_{H(P)}: H(P)\to \Phi_P\left(H(P)\right)$ is bijective.   
Hence, we have the following:  
\begin{lemma}\label{lemma11}
For any $P\in S^n$, the mapping 
\[
\Phi_P|_{H(P)}: H(P)\to \Phi_P\left(H(P)\right)
\]
is a $C^\infty$ diffeomorphism.   
\end{lemma}
\par 
Since $\Psi_P(\mathbf{x})\cdot P>0$ for any 
$\mathbf{x}\in S^n-\{\pm P\}$, by Property \ref{property1},  
it follows that 
$ped_{\gamma,P}(I)$ $\subset$ ${H}(P)$. 
Therefore, by Lemma \ref{lemma11}, 
the assertion (1) of Theorem \ref{thm1} holds. 
\hfill $\Box$
\begin{note}\label{note2}
The proof given above actually shows that if $P$ is located outside 
${\bf u}_n(I)$, then $ort_{\gamma, P}$ and $ped_{\gamma, P}$ are 
$\mathcal{L}$-equivalent.   
Hence, the assertion (1) of Theorem \ref{thm1} 
can be improved to the global version as follows.   
\end{note}
\begin{theorem}\label{thm2}
Let $\gamma: I\to S^n$ be a spherical unit speed curve and let 
$P$ be a point of $S^n$.    
 Then, 
The orthotomic curve $ort_{\gamma, P}: I\to S^n$ is
$\mathcal{L}$-equivalent to
 $ped_{\gamma, P}: I\to S^n$ if and only if 
$P\not\in {\bf u}_n(I)$.
\end{theorem}
\subsection{Proof of the assertion (2) of Theorem \ref{thm1}}
Let $P$ be the point satisfying 
$P={\bf u}_n(s_0)$ or $P=-{\bf u}_n(s_0)$ and 
let $U$ be a sufficiently small neighborhood of 
$P$ in $S^n$.     
By the following lemma, 
the equality given in Lemma \ref{lemma8} can be extended to the case 
that $ped_{\gamma, P}$ can not be defined.         
\begin{lemma}\label{lemma composition}
For any ${\bf x}\in U-\{P\}$, the following holds:   
\[
\Phi_P\circ\Psi_P({\bf x})=-\Phi_P({\bf x}).   
\]
\end{lemma}    
\proof 
\begin{eqnarray*}
\Phi_P\circ \Psi_P({\bf x}) 
& = & 
2\left(P\cdot \Psi_P({\bf x})\right)\Psi_P({\bf x})-P \\ 
{ } & = & 
2\frac{1}{\sqrt{1-(P\cdot {\bf x})^2}}
\left(1-(P\cdot {\bf x})^2\right)\Psi_P({\bf x})-P \\ 
{ } & = & 
2\sqrt{1-(P\cdot {\bf x})^2}\frac{1}{\sqrt{1-(P\cdot {\bf x})^2}}
\left(P-(P\cdot {\bf x}){\bf x}\right)-P \\ 
{ } & = & 
2\left(P-(P\cdot {\bf x}){\bf x}\right)-P \\ 
{ } & = & 
P-2(P\cdot {\bf x}){\bf x} \\ 
{ } & = & 
-\Phi_P({\bf x}).   
\end{eqnarray*}
\hfill $\Box$
\par 
Since the assertion is concerned with map-germ at $s_0$, we may assume 
$I$ is sufficiently small so that $I\subset \pm{\bf u}_n^{-1}(U)$.     
By Lemma \ref{lemma composition}, if $\pm{\bf u}_n(s)$ is inside $U-\{P\}$, 
then the following holds:   
\[
ort_{\gamma, P}(s)=\Phi_P\circ ped_{\gamma, P}(s) 
= \Phi_P\circ \Psi_P\circ {\bf u}_n(s) 
= -\Phi_P\circ {\bf u}_n(s).    
\] 
Since all of $ort_{\gamma, P}: I\to S^n$, ${\bf u}_n: I\to S^n$ 
and $\Phi_P: S^n\to S^n$ are continuous mappings, 
the above equality 
$ort_{\gamma, P}(s)=-\Phi_P\circ {\bf u}_n(s)$ can be extended 
to the $s$ satisfying $P=\pm{\bf u}_n(s)$.  
In other words, two mappings  
$ort_{\gamma, P}: I\to S^n$ and 
$-\Phi_P\circ {\bf u}_n: I\to S^n$ are exactly the same mappings.   
Moreover, since $U$ is 
sufficiently small and any point of $U$ is a regular point of 
$\Phi_P$, it follows that $-\Phi_P|_U: U\to -\Phi_P(U)$ is a 
$C^\infty$ diffeomorphism.    
Therefore, $ort_{\gamma, P}$ and ${\bf u}_n$ are $\mathcal{L}$-equivalent.   
\hfill 
$\Box$
\begin{note}\label{note3}
Consider a unit-speed curve $\gamma: (-2, 2)\to S^n$ such that 
${\bf u}_n$ is an embedding and ${\bf u}_n(-1)=-{\bf u}_n(1)$ is 
satisfied.   Set $P={\bf u}_n(1)$.   Then, 
$ort_{\gamma, P}(-1)=ort_{\gamma, P}(1)=-P$.    
Hence, $ort_{\gamma, P}: (-2, 2)\to S^n$ is not injective.   
Therefore, $ort_{\gamma, P}$ and ${\bf u}_n$ are 
not $\mathcal{A}$-equivalent.      
This example shows that the assertion (2) of Theorem \ref{thm1} 
cannot be improved to the global version.      
\end{note}
\section{Appendix: Plane pedal curves and plane orthotomic curves}
\label{section 5}
\begin{figure}[h]
\centering
\begin{tikzpicture}
  \draw [<-] (1,1)node[below]{$T(s)$} -- (5,5);
  \draw [thick] (1.5,0.5) to [out=75,in=225] (2.5,2.5);
  \draw [thick](2.5,2.5) to [out=45,in=190] (4.5,3.5)node[below]{$\gamma$};
  \draw [->] (2.5,2.5)--(4,1)node[below]{$N(s)$};
  \draw [fill] (5,2) circle [radius=0.05] node[below]{$P$};
  \draw [fill] (3.5,3.5) circle [radius=0.05] node[above right]{\quad
$ped_{\gamma,P}(s)=\left((\gamma(s)-P) \cdot N(s)\right)N(s)+P$};
  \draw [fill] (2,5) circle [radius=0.05];
  \draw [dashed] (5,2)--(3.5,3.5);
  \draw [dashed] (3.5,3.5)--(2,5) node[above]
{$ort_{\gamma,P}(s)=2\left((\gamma(s)-P) \cdot N(s)\right)N(s)+P$};
  \draw [-] (3.6,3.6)--(3.7,3.5);
  \draw [-] (3.6,3.4)--(3.7,3.5);
  \draw [-] (2.6,2.6)--(2.7,2.5);
  \draw [-] (2.6,2.4)--(2.7,2.5);
\end{tikzpicture}
\centering
\caption{$ped_{\gamma, P}(s)$ and $ort_{\gamma, P}(s)$ in $\mathbb{R}^2$.}
\label{figure1}
\end{figure}
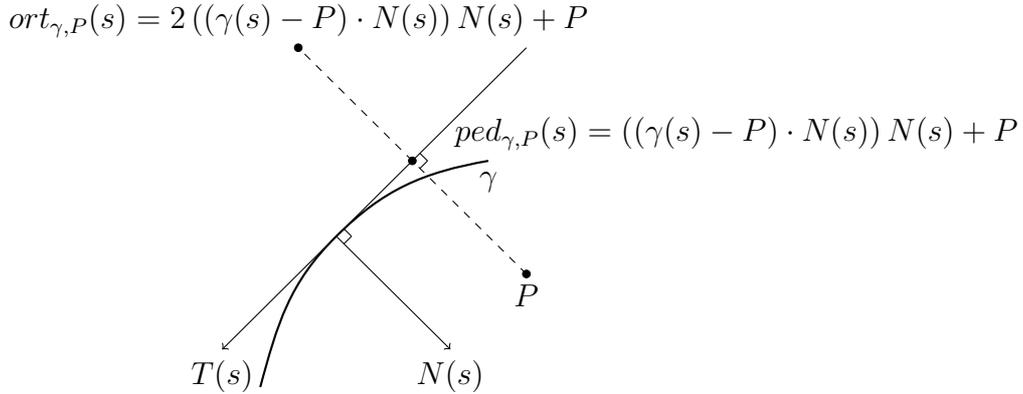

Let $\gamma: I\to \mathbb{R}^2$ be a unit-speed 
and let $P$ be a point of  $\mathbb{R}^2$.  
Following J.W.~Bruce and P.J.~Giblin \cite{brucegiblin}, 
the pedal curve of $\gamma$ relative to $P$ and the orthotomic curve 
of $\gamma$ relative to $P$ are reviewed.   
\begin{definition}[\cite{brucegiblin}]\label{planecurvecase}
{\rm
\begin{enumerate}
\item[(1)]   The {\it pedal curve of }$\gamma$ relative to $P$, 
denoted by $ped_{\gamma, P}: I\to \mathbb{R}^2$, is defined as follows:   
\[
ped_{\gamma, P}(s)-P=\left(\left(\gamma(s)-P\right)\cdot N(s)\right)N(s).
\]  
\item[(2)]   The {\it orthotomic curve of }$\gamma$ relative to $P$, 
denoted by $ort_{\gamma, P}: I\to \mathbb{R}^2$, is defined as follows:   
\[
ort_{\gamma, P}(s)-P=2\left(\left(\gamma(s)-P\right)\cdot N(s)\right)N(s).
\]  
\end{enumerate} 
Here, $N(s)$ stands for the unit normal vector to $\gamma$ at $\gamma(s)$.
}
\end{definition} 
%
Let $L: \mathbb{R}^2\to \mathbb{R}^2$ be the mapping defined by 
$L({\bf x})=2{\bf x}-P$.   Then, $L$ is clearly a $C^\infty$ diffeomorphism and 
$ort_{\gamma, P}=L\circ ped_{\gamma,P}$ (see Figure \ref{figure1}).     
Hence, $ort_{\gamma, P}$ is always $\mathcal{L}$-equivalent to 
$ped_{\gamma, P}$.  
\section*{Acknowledgement}
The second author was 
supported
by JSPS KAKENHI Grant Number 17K05245.

\end{document}